\documentclass[12pt]{article}
\usepackage{amssymb}
\usepackage{amsthm}
\usepackage{amsmath}
\usepackage{graphicx}
 \newtheorem{thm}{Theorem}[section]
 \newtheorem{cor}[thm]{Corollary}
 \newtheorem{lem}[thm]{Lemma}
 
 \theoremstyle{definition}
 
 \newtheorem{rem}[thm]{Remark}
 \numberwithin{equation}{section}

\begin{document}
\title{A Liouville Theorem for the Axially-symmetric Navier-Stokes Equations}
\author{Zhen Lei\footnote{School of
Mathematical Sciences; LMNS and Shanghai Key
  Laboratory for Contemporary Applied Mathematics, Fudan
  University, Shanghai 200433, P. R. China. {\it Email:
 leizhn@yahoo.com, zlei@fudan.edu.cn}}\and
 Qi S. Zhang\footnote{Department of Mathematics, University of
 California, Riverside, CA 92521, USA. {\it Email: qizhang@math.ucr.edu}}}
\date{\today}
\maketitle

\begin{abstract}
Let $v(x, t)= v^r e_r + v^\theta e_\theta + v^z e_z$ be a solution
to the three-dimensional incompressible axially-symmetric
Navier-Stokes equations. Denote by $b = v^r e_r + v^z e_z$ the
radial-axial vector field. Under a general scaling invariant
condition on $b$, we prove that the quantity $\Gamma = r v^\theta$
is H\"older continuous at $r = 0$, $t = 0$. As an application, we
give a partial proof of a conjecture on Liouville property by
Koch-Nadirashvili-Seregin-Sverak in \cite{KNSS} and Seregin-Sverak
in \cite{SS}. As another application, we prove that if $b \in
L^\infty([0, T], BMO^{-1})$, then $v$ is regular. This provides an
answer to an open question raised by Koch and Tataru in
\cite{KochTataru} about the uniqueness and regularity of
Navier-Stokes equations in the axially-symmetric case.

\end{abstract}


\section{Introduction}

In this paper we study the three-dimensional incompressible
axially-symmetric Navier-Stokes equations. In cylindrical
coordinates, the velocity field $v = v(x, t)$ is of the form
\begin{equation}\nonumber
v(x, t) = v^r(r, z, t)e_r +  v^\theta(r, z, t)e_\theta +  v^z(r,
z, t)e_z.
\end{equation}
Here and throughout the paper, we write $x = (x_1, x_2, z)$, $r =
r(x) = \sqrt{x_1^2 + x_2^2}$ and
\begin{equation}\nonumber
e_r = e_r(x) = \begin{pmatrix}\frac{x_1}{r}\\ \frac{x_2}{r}\\ 0
  \end{pmatrix},
\quad e_\theta = e_\theta(x) = \begin{pmatrix}- \frac{x_2}{r}\\ \frac{x_1}{r}\\
0
  \end{pmatrix},
\quad e_z = e_z(x) = \begin{pmatrix}0\\ 0\\ 1
  \end{pmatrix}
\end{equation}
are the three orthogonal unit vectors along the radial, the
angular, and the axial directions respectively. The radial,
angular (or swirl) and axial components $v^r$, $v^\theta$ and
$v^z$ of the velocity field are governed by (see, for instance,
\cite{MajdaB})
\begin{equation}\label{axi-NS}
\begin{cases}
\partial_tv^r + b\cdot\nabla v^r - \frac{(v^\theta)^2}{r}
  + \partial_rp = \big(\Delta - \frac{1}{r^2}\big)v^r,\\[-4mm]\\
\partial_tv^\theta + b\cdot\nabla v^\theta + \frac{v^rv^\theta}{r}
  = \big(\Delta - \frac{1}{r^2}\big)v^\theta,\\[-4mm]\\
\partial_tv^z + b\cdot\nabla v^z + \partial_zp = \Delta v^z,\\[-4mm]\\
b = v^re_r + v^ze_z,\quad \nabla\cdot b = \partial_rv^r +
\frac{v^r }{r} + \partial_zv^z = 0.
\end{cases}
\end{equation}
Here without loss of generality, we have set the viscosity
constant to be unit.

A special feature of the axially-symmetric Navier-Stokes equations
is that  the quantity $\Gamma = r v^\theta (x, t)$ satisfies an
parabolic equation with singular drift terms:
\begin{equation}\label{1.5}
\partial_t\Gamma + b\cdot\nabla\Gamma + \frac{2}{r}\partial_r\Gamma = \Delta\Gamma.
\end{equation}
We remark that $\Gamma$ enjoys the maximal principle. For this
reason the axially-symmetric case appears more tractable than the
full three-dimensional problem.

Nevertheless, it is well-known that global regularity of the
three-dimensional incompressible Navier-Stokes equations is still
wide open even in the axially-symmetric case. But if the swirl
component of the velocity field $v^\theta$ is trivial,
independently, Ladyzhenskaya \cite{L} and Uchoviskii \& Yudovich
\cite{UY} proved that weak solutions are regular for all time (see
also \cite{LMNP}). Recently, tremendous efforts and interesting
progresses have been made on the regularity problem of the
axially-symmetric Navier-Stokes equations with a general
non-trivial swirl. For example, in \cite{CSYT2008, CSYT2009},
Chen-Strain-Tsai-Yau proved, among other things, that the suitable
weak solutions are smooth if the velocity field $v$ satisfies
$r|v| \leq C_\ast < \infty$. Their method is based on the
classical results by Nash \cite{Nash}, Moser \cite{Moser} and De
Giorgi \cite{DG}. In \cite{KNSS}, Koch-Nadirashvili-Seregin-Sverak
proved the same result using Liouville type theorem for ancient
solutions of Navier-Stokes equations. See also \cite{SS} for a
local version.


A velocity field is called an ancient solution if it exists in the
time interval $(-\infty, 0]$, and it satisfies the Navier-Stokes
equation in certain sense. A well known fact is that ancient
solutions represent structures of singularity of evolution
equations, which makes the study of ancient solutions an important
topic.

In this paper, we study the axially-symmetric Navier-Stokes
equations
 under a more general assumption on the radial-axial velocity vector
$b$. To be precise, we consider $b$ such that
\begin{equation}\label{1.1}
b = b_1 + b_2 + b_3,\quad \nabla\cdot b_1 = \nabla\cdot b_2 =
\nabla\cdot b_3 = 0,
\end{equation}
where
\begin{equation}\label{1.2}
{\rm HSE}(b_1) \leq C_\ast,\quad b_2 = \nabla\times B,\ \sup_{-T <
t < 0}\|B\|_{{\rm BMO }} \leq C_\ast,\quad  \sup_{-T < t < 0, x
\in \mathbb{R}^3}r|b_3| \leq C_\ast.
\end{equation}
Some motivation and explanation for the condition and notations
are in order. Here $[-T, 0]$ is the time interval where a solution
exists. We often take $T=1$ for convenience. The number $C_\ast$
is an arbitrary positive constant and  ${\rm HSE}(b_1)$ is called
"the hollowed scaled energy", defined by
\begin{equation}\label{1.3}
{\rm HSE}(b_1) = \sup_{0 < R < R_0}\dot{E}_R(b_1),\quad
\dot{E}_R(b_1) = \sup_{- R^2 < t <
0}\frac{1}{R}\int_{B_{2R}/B_{R/8}}|b_1(\cdot, t)|^2dx.
\end{equation}  Here $R_0$ is a positive number often taken as
$1$.

 We
use
\begin{equation}
\label{Enorm} \Vert b \Vert_E = {\rm HSE}(b_1) + \sup_{-T < t < 0}
\|B\|_{{\rm BMO}} + \sup_{-T < t < 0} r |b_3|
\end{equation} to denote the controlling quantity of $b$
throughout the paper. Here $[-T, 0]$ is the time interval of
concern, which may be shifted or scaled.  The linear space
consisted of those $b$ such that $\Vert b \Vert_E<\infty$ is
called space $E$. The results in this paper depends on $b$ only in
terms of $\Vert b \Vert_E$. We will use a positive function
$K(\Vert b \Vert_E)$ to denote such a dependence, whose precise
value may change from line to line. Notice that the space $E$
contains $BMO^{-1}$ which is the largest known space in which the
Navier-Stokes equations are well-posed. See the interesting work
by Koch and Tataru \cite{KochTataru}. Another feature is that the
condition on $b_1$ is imposed only on some subdomain of the space
time cube. Outside of the subdomain, there is no restriction on
$b_1$. With a little bit more efforts, we can also just impose
conditions on part of the space time for $b_2$ and $b_3$ too. But
here we do not pursue that.

Our first result states that $\Gamma =r v^\theta$ is H\"older
continuous at $r = 0, t=0$ if the radial-axial velocity field $b$
satisfies \eqref{1.1} and \eqref{1.2}. The H\"older continuity
depends on $b$ only through $\| b \|_E$.

\begin{thm}\label{Holder}
Given a number $L>0$, let $v=v(x, t)$, $(x, t) \in Q_L \equiv B(x,
L) \times [-L^2, 0]$ be a $L^\infty_{{\rm loc}}(Q_L)$ weak
solution to the three-dimensional axially-symmetric Navier-Stokes
equations \eqref{axi-NS}. Suppose that the radial-axial velocity
field $b$ satisfies \eqref{1.1}-\eqref{1.2}. Then $\Gamma =
rv^\theta$ is H${\rm \ddot{o}}$lder continuous at $(0, 0)$
uniformly. i.e. There exist positive constants $\alpha$ and $C$,
depending only on $\| b \|_E$, such that, for all $(x, t) \in
Q_{L/2}$, it holds
\[
|\Gamma(x,t) -\Gamma(0, 0)| \le C [(|x|+\sqrt{|t|})/L]^{\alpha} \,
\sup_{Q_L} |\Gamma|.
\]
\end{thm}

Our proof is inspired by  \cite{CSYT2008} where the authors had
proved a version of the above theorem under the assumption $r|v|
\leq C_\ast$ using an De Giorgi type argument (see also
\cite{CSYT2009} for the method based on the direct estimation of
an evolution kernel). Here we will treat the more general $b$
using a Nash type method in a uniform way. We will first establish
a local maximum estimate for solutions of \eqref{1.5} in terms of
the controlling constant $C_\ast$ for $b$ in (\ref{1.2}). This is
done by using Moser's iteration method and De Giorgi type energy
estimate, exploiting the structure of $b$. Similar argument has
appeared in Zhang \cite{Zhang}  and Chen-Strain-Tsai-Yau
\cite{CSYT2008, CSYT2009} where $b$ is some form  bounded function
or $|r b(x, t)| \le C_\ast$.  Then we apply the Nash type method
to prove the H\"older continuity of $\Gamma$. One handy tool which
allows to treat more general type of vector fields $b$ is a simple
two dimensional integration by parts argument (\ref{imbed}).
Another tool is the John-Nirenberg inequality for BMO functions,
which was first employed by Friedlander and  Vicol \cite{FV}, and
 also by Seregin, Silvestre, Sverak,
Zlatos \cite{SSSA} to treat the linear heat equation with $\Delta
u + b \nabla u -\partial_t u=0$ with $b \in L^\infty([0, T],
BMO^{-1})$. They prove H\"older continuity of weak solutions to
this equation. We also utilize the role played the stream
function, which helps to do integration by parts one more time.
Let $v$ be a velocity field. We recall that a function $B$ is
called a stream function of $v$ if $v=\nabla \times B$.

The main significance of Theorem \ref{Holder} is that it deduces
the next two theorems.  One of them gives a partial answer to an
open question in \cite{KNSS} on Liouville properties. The other
one establishes a condition on $b$ such that solutions to
axially-symmetric Navier-Stokes equations are regular. This
regularity condition does not involve Lebesgue integral on $b$ or
absolute value of $b$, which may allow the capturing of more
oscillatory functions.

\begin{thm}\label{Liouville}
Let $v=v(x, t)$ be a bounded, weak ancient solution to
(\ref{axi-NS}). Suppose also $r |v^\theta|$ is bounded and the
stream function is a BMO function. Then $v \equiv 0$.
\end{thm}

\begin{rem} The authors of \cite{KNSS} stated a conjecture on
Liouville type theorem for the axially-symmetric Navier-Stokes
equations: bounded, mild, ancient velocity fields are constants.
The authors in \cite{KNSS} proved such kind of Liouville theorems
in the three-dimensional axially-symmetric case without swirl, or
under the condition $r|v|$ being bounded. The above theorem, under
the conditions that $r |v^\theta|$ is bounded and the stream
function is a BMO function, gives a proof of this conjecture.

Recall that $r v^\theta$ is scaling invariant and it also
satisfies the maximum principle. Therefore its boundedness is a
mild restriction. A bounded function is obviously a BMO function.
Although a bounded velocity field may not have a bounded stream
function in general, a boundedness assumption on the stream
function is also very mild since one expects it to hold in most
natural cases when the velocity is bounded.
\end{rem}

In \cite{KochTataru} p25, Koch and Tataru wrote that there has
been a strong interest in obtaining well-posedness of
Navier-Stokes equations assuming a ${\rm BMO^{-1}}$ space
condition. They raised the question of uniqueness and regularity
for solutions in the space $L^\infty( [0, T], BMO^{-1})$. They
also proved uniqueness for small solutions in such a space. Our
Theorem \ref{thregcon} answers this question in the
axially-symmetric case. We should mention that decay and
analyticity estimates for solutions to the Navier-Stokes equations
evolving from small initial data in $BMO^{-1}$ are studied
independently by Miura and Sawada \cite{Mih}  and Germain,
Pavlovi${\rm \acute{c}}$ and Staffilani \cite{GPS}.

\begin{thm}\label{thregcon}
Let $v=v(x, t)$ be a suitable weak solution to (\ref{axi-NS}) in
the space time region ${\bf R}^3 \times [0, T]$. Assume  that the
initial value satisfies $v(\cdot, 0) \in L^2$ , $|r v^\theta(x,
0)|<C$. Suppose also $v(\cdot, t) = \nabla \times B(\cdot, t)$
with $\sup_{0 < t < T}\|B(\cdot, t)\|_{{\rm BMO}} \leq C_\ast$.
Then $v$ is smooth in ${\bf R}^3 \times (0, T]$.  Here $C$ and
$C_\ast$ are arbitrary positive constants.
\end{thm}

\begin{rem}
Note condition  $|r v^\theta(x, 0)|<C$ is only on the initial
value. It can also be dropped by a approximation argument. We will
not seek the full generality this time.
\end{rem}

\begin{rem}
In \cite{ESS}, Escauriaza, Seregin and Sverak proved that
$L^\infty_TL^3(Q)$ solutions to the Navier-Stokes equations are
regular, which is the highly non-trivial borderline case of
Serrin's criterion. Their proof is based on the method of backward
uniqueness and unique continuation together with a blowup
argument.  Since $L^3$ is imbedded into ${\rm BMO^{-1}}$, our
Theorem \ref{thregcon} also provides a new and simpler proof to
such a criterion in the axially-symmetric case.
\end{rem}

Before ending the introduction, let us mention some other related
results on axially-symmetric Navier-Stokes equations. In the
presence of swirl, there is the paper by J. Neustupa \& M. Pokorny
\cite{NP}, proving the regularity of one component (either $v^r$
or $v^{\theta}$) implies regularity of the other components of the
solution. Also proving regularity is the work of Q. Jiu \& Z. Xin
\cite{JX} under an assumption of sufficiently small zero dimension
scaled norms. We would also like to mention the regularity results
of D. Chae \& J. Lee \cite{CL} who prove regularity results
assuming finiteness of another zero dimensional integral. On the
other hand, G. Tian \& Z. Xin \cite{TX} constructed a family of
singular axially symmetric solutions with singular initial data;
T. Hou \& C. Li \cite{HL} found a special class of global smooth
solutions. See also a recent extension: T. Hou, Z. Lei \& C. Li
\cite{HLL}.

The paper is organized as follows: In section 2 we establish a
local maximum estimate using De Giorgi type energy method and
Moser's iteration method. Based on the local maximal estimate, we
obtain the H${\rm \ddot{o}}$lder continuity of $\Gamma$ and prove
Theorem \ref{Holder} by Nash's method in section 3. The argument
is based on \cite{CSYT2008, CSYT2009}. Then in section 4 we prove
our Theorem \ref{Liouville} and Theorem \ref{thregcon}, using
Theorem \ref{Holder} and some new blow up arguments. The main idea
is that a possible singularity falls only into two types. Type I
singularity can be scaled into an axially-symmetric, bounded,
ancient mild solution. Type II can be scaled to a two dimensional
ancient solution. Then we show that either type leads to a
contradiction with the assumption that the stream function is in
the BMO space. In the process the two dimensional Liouville
theorem in \cite{KNSS} plays an important role.

\section{Local Maximum Estimate}

In this section we prove a local maximum estimate of $\Gamma$
using Moser's iteration method in proving the parabolic Harnack's
inequality. These estimates will be used to obtain H${\rm
\ddot{o}}$lder continuity of $\Gamma$ in next section. The main
idea is to exploit the divergence-free property of $b(x, t)$ and
to construct a special cut-off function. We also learned from
\cite{CSYT2008, CSYT2009} where the authors treated the term
$\frac{2}{r} \partial_r \Gamma$ in the equation for $\Gamma$.

We first derive an energy estimate of  De Giorgi type for
\eqref{1.5}. For this purpose we need a refined cut-off function.
Set $\frac{1}{2}\leq
\sigma_2 < \sigma_1 \leq 1$ and choose $\psi(y, s) =
\phi(|y|)\eta(s)$ to be a smooth cut-off function satisfying:
\begin{equation}\label{E2-1}
\begin{cases}
{\rm supp}\phi \subset B(\sigma_1),\ \phi = 1\
  {\rm on}\ B(\sigma_2),\ 0 \leq \phi \leq 1,\\[-4mm]\\
{\rm supp} \eta \subset \big( - (\sigma_1)^2, 0\big],\ \eta(s)=1\
{\rm
  on}\ \big( - (\sigma_2)^2, 0\big],\ 0\leq \eta \leq 1,\\[-4mm]\\
|\eta '| \lesssim \frac{1}{(\sigma_1 - \sigma_2)^2},\
  \big|\frac{\nabla\phi}{\sqrt{\phi}}\big| \lesssim \frac{1}{\sigma_1 - \sigma_2},\quad
  \big|\nabla\big(\frac{\nabla\phi}{\sqrt{\phi}}\big)\big| \lesssim \frac{1}{(\sigma_1 - \sigma_2)^2}.
\end{cases}
\end{equation}
Here as usual we use $A \lesssim B$ to denote the inequality $A
\leq CB$ for an absolute positive constant $C$. Such a cut-off
function $\phi$ can be simply chosen as a square of a standard
cut-off function. We will also use the following notations for
domains. Let $R>0$, we write $B_R=B(0, R)$ and
\begin{eqnarray}\nonumber
P(R) = B_R \times (- R^2, 0], \quad P(R_1, R_2) = B_{R_1}/ B_{R_2}
\times (- R_1^2, 0]\ {\rm for}\ R_1 > R_2.
\end{eqnarray}

Consider the functions $f = |\Gamma|^q$ , $q
> \frac{1}{2}$ and the cut-off functions
$\psi_R(y, s) = \phi_R(y)\eta_R(s) =
\phi(\frac{y}{R})\eta\big(\frac{s}{R^2}\big)$. Testing \eqref{1.5}
by $q|\Gamma|^{2q - 2}\Gamma\psi_R^2$ gives
\begin{eqnarray}\label{E2-2}
\frac{1}{2}\iint
  \Big(\partial_sf^2 + (b\cdot\nabla)f^2 + \frac{2}{r}\partial_rf^2\Big)\psi_R^2dyds = q\iint
  \Delta\Gamma|\Gamma|^{2q - 2}\Gamma\psi_R^2dyds.
\end{eqnarray}
Using Cauchy-Schwarz's inequality and integration by parts, we
compute that
\begin{eqnarray}\nonumber
&&q\iint(\Delta
  \Gamma)|\Gamma|^{2q - 2}\Gamma\psi_R^2dyds = q\iint(\Delta
  |\Gamma|)|\Gamma|^{2q - 1}\psi_R^2dyds\\\nonumber
&&= - q\iint\Big((2q - 1)|\nabla\Gamma|^2\Gamma^{2q -
  2}\psi_R^2 + \nabla\psi_R^2\cdot|\Gamma|^{2q - 1}\nabla
  |\Gamma|\Big)dyds\\\nonumber
&&= - \iint\Big((2 - \frac{1}{q}
  )|\nabla f|^2\psi_R^2 + 2\psi_R\nabla\psi_R\cdot f\nabla
  f\Big)dyds\\\nonumber
&&= - \iint\Big((2 - \frac{1}{q}
  )|\nabla f|^2\psi_R^2 + 2f\nabla\psi_R\cdot \nabla
  (f\psi_R) - 2f^2|\nabla\psi_R|^2\Big)dyds\\\nonumber
&&\lesssim - \iint
  |\nabla(f\psi_R)|^2dyds + \iint f^2|\nabla\psi_R|^2dyds
\end{eqnarray}
and
\begin{eqnarray}\nonumber
&&\frac{1}{2}\iint\psi_R^2
  \partial_sf^2dyds\\\nonumber
&&= \frac{1}{2}\int_{B(\sigma_1R)}\psi_R^2
  f^2(\cdot, t)dy - \frac{1}{2}\iint f^2\partial_s\psi_R^2dyds.
\end{eqnarray}
Moreover, by the fact that $\Gamma = 0$ on the axis $r = 0$, we
have
\begin{eqnarray}\nonumber
\iint\frac{1}{r}\partial_rf^2\psi_R^2dyds =
\iint\partial_rf^2\psi_R^2drdzd\theta ds = \iint
f^2\partial_r\psi_R^2drdzd\theta ds.
\end{eqnarray}
Consequently, using \eqref{E2-1}, we have
\begin{eqnarray}\label{E2-3}
&&\frac{1}{2}\int\psi_R^2
  f^2(\cdot, t)dy + \iint|\nabla(f\psi_R)|^2dyds\\\nonumber
&&\lesssim \frac{1}{(\sigma_1 - \sigma_2)^2R^2}
  \iint_{P(\sigma_1R)}f^2dyds  - \frac{1}{2}\iint
  (b\cdot\nabla f^2)\psi_R^2dyds.
\end{eqnarray}

Now we start to treat the drift term involving $b=b_1+b_2+b_3$.
For $R_1
> R_2$, let us denote that
\begin{equation}\nonumber
\dot{E}(R_1, R_2, b) = \sup_{- R_1^2 \leq t \leq 0}\frac{1}{R_1 -
R_2}\int_{B_{R_1}\setminus B_{R_2}}|b(\cdot, t)|^2dx.
\end{equation}
By \eqref{E2-1} and the divergence-free properties of the velocity
field $b_1(x, t)$, we have
\begin{eqnarray}\nonumber
&&- \frac{1}{2}\iint(b_1\cdot\nabla f^2)\psi_R^2dyds =
  \iint_{P(\sigma_1R,\ \sigma_2R)}
  b_1\cdot\frac{\nabla\phi_R}{\phi_R^{\frac{1}{2}}}(\psi_Rf)^{\frac{3}{2}}
  (\eta_Rf)^{\frac{1}{2}}dyds\\\nonumber
&&\lesssim \frac{1}{(\sigma_1 - \sigma_2)R}\int\|b_1\|_{L^2
  \big(B(\sigma_1R,\ \sigma_2R)\big)}\|\psi_Rf\|_{L^6\big(B(\sigma_1R)
  \big)}^{\frac{3}{2}}\|f\|_{L^2\big(B(\sigma_1R)\big)}^{\frac{1}{2}}ds\\\nonumber
&&\lesssim \Big(\frac{\dot{E}(\sigma_1R, \sigma_2R, b_1)}
  {(\sigma_1 - \sigma_2)R}\Big)^{\frac{1}{2}}\|\psi_Rf\|_{L^2_tL^6\big( P(\sigma_1R)
  \big)}^{\frac{3}{2}}\|f\|_{L_t^2L^2\big(P(\sigma_1R)\big)}^{\frac{1}{2}}.\\\nonumber
\end{eqnarray} Therefore
\begin{eqnarray}
\label{b1=1} &&- \frac{1}{2}\iint(b_1\cdot\nabla
f^2)\psi_R^2dyds\\\nonumber
 &&\lesssim
\frac{\dot{E}(\sigma_1R, \sigma_2R, b_1)^2}{(\sigma_1 -
  \sigma_2)^2R^2}\|f\|_{L^2_tL^2\big(P(\sigma_1R)\big)}^2 +
  \frac{1}{8}\iint_{P(\sigma_1R)}|\nabla(\psi_Rf)|^2dyds.
\end{eqnarray}

Next we treat the term involving $b_2$.  Let $\bar{B}=\bar{B}(t)$
be the average of $B(\cdot, t)$ in $B_R$.  Then
\begin{eqnarray}\nonumber
&&- \frac{1}{2}\iint(b_2\cdot\nabla f^2)\psi_R^2dyds\\\nonumber
&&= \iint_{P(\sigma_1R,\ \sigma_2R)}
  (B-\bar{B}(t)) \cdot\nabla\times\big(\frac{\nabla\phi_R}{\phi_R^{\frac{1}{2}}}(\psi_Rf)^{\frac{3}{2}}
  (\eta_Rf)^{\frac{1}{2}}\big)dyds\\\nonumber
&&\lesssim \big\|\nabla\big(\frac{\nabla\phi_R}
  {\sqrt{\phi_R}}\big)\big\|_{L^\infty_tL^2\big(P(\sigma_1R,\ \sigma_2R)\big)}
  \|\psi_Rf\|_{L^2_tL^6\big( P(\sigma_1R)\big)}^{\frac{3}{2}}\|(B-\bar{B}) f\|_{L^2_tL^2
  \big(P(\sigma_1R)\big)}^{\frac{1}{2}}\\\nonumber
&&\quad +\ \iint_{P(\sigma_1R,\ \sigma_2R)}
 (B-\bar{B}) \cdot\big[\frac{\nabla\phi_R}{\phi_R^{\frac{1}{2}}}\times\nabla
  \frac{(\psi_Rf)^2}{\sqrt{\phi_R}}\big]dyds.
\end{eqnarray} Therefore
\begin{eqnarray}\nonumber
  \\\nonumber
  &&- \frac{1}{2}\iint(b_2\cdot\nabla f^2)\psi_R^2dyds\\\nonumber
&&\lesssim \Big(\frac{1}{(\sigma_1 - \sigma_2)R}
  \|(B-\bar{B}) f\|_{L^2_tL^2\big(P(\sigma_1R)\big)}\Big)^{\frac{1}{2}}
  \|\nabla(\psi_Rf)\|_{L^2_tL^2\big( P(\sigma_1R)
  \big)}^{\frac{3}{2}}\\\nonumber
&&\quad +\ \frac{1}{(\sigma_1 - \sigma_2)R}
  \|(B-\bar{B}) f\|_{L^2_tL^2\big(P(\sigma_1R)\big)}\|\nabla(\psi_Rf)\|_{L^2_tL^2\big( P(\sigma_1R)
  \big)}\\\nonumber
&&\quad +\ \Big(\frac{1}{(\sigma_1 - \sigma_2)R}
  \|(B-\bar{B}) f\|_{L^2_tL^2\big(P(\sigma_1R)\big)}\Big)^2.\\\nonumber
\end{eqnarray}
Hence
\begin{equation}\label{b2=1}
\begin{aligned}
 - \frac{1}{2}&\iint(b_2\cdot\nabla f^2)\psi_R^2dyds\lesssim \frac{1}{8}\iint_{P(\sigma_1R)}|\nabla(\psi_Rf)|^2dyds\\
& +\ \frac{1}{(\sigma_1 -
  \sigma_2)^2R^2}\|(B-\bar{B}) f\|_{L^2_tL^2\big(P(\sigma_1R)\big)}^2.
\end{aligned}
\end{equation} To control the last expression, we need to recall
the well-known John-Nirenberg inequality  for BMO functions (see
\cite{JohnN} or \cite{Stein}): for any $p \in (0, \infty)$,
\begin{equation}
\label{BMOLp} \Vert B(\cdot, t) - \bar{B}(t) \Vert_{L^p(B_R)} \le
C_p \Vert B(\cdot, t) \Vert_{BMO} |B_R|^{1/p}.
\end{equation}
Taking $p=6$ in the above inequality,  we have
\[
\begin{aligned}
&\|(B-\bar{B}) f\|_{L^2_tL^2\big(P(\sigma_1R)\big)}\le \|
f\|_{L^3_tL^3\big(P(\sigma_1R)\big)} \| B-\bar{B}
\|_{L^6_tL^6\big(P(\sigma_1R)\big)} \\
&\le  \| f\|_{L^3_tL^3\big(P(\sigma_1R)\big)} \Vert B
\Vert_{L^\infty_t BMO} |B_R|^{\frac{1}{6}} R^{\frac{1}{3}}.
\end{aligned}
\]Plug this into (\ref{b2=1}), we deduce
\begin{equation}
\label{b2=2}
\begin{aligned}
 - \frac{1}{2}&\iint(b_2\cdot\nabla f^2)\psi_R^2dyds\\
&\lesssim \frac{ \Vert B \Vert^2_{L^\infty_t BMO}
R^{\frac{5}{3}}}{(\sigma_1 -
  \sigma_2)^2R^2} \| f\|^2_{L^3_tL^3\big(P(\sigma_1R)\big)} +
  \frac{1}{8}\iint_{P(\sigma_1R)}|\nabla(\psi_Rf)|^2dyds.
\end{aligned}
\end{equation}

The term involving $b_3$ has been treated in \cite{CSYT2009}. Here
we give an alternative proof for completeness and simplicity.
\begin{equation}
\label{imbed}
\begin{aligned}
&| \frac{1}{2}\iint(b_3 \cdot\nabla f^2)\psi_R^2dyds| \lesssim |
\frac{1}{2}\iint(b_3 \cdot \nabla (\psi_R^2) f^2) dyds|\\
&= |\iint(b_3 \cdot \nabla (\psi_R) \psi_R f^2) dyds| \lesssim
|\iint(\frac{1}{r} |\nabla \psi_R| \psi_R f^2) r drd\theta ds|\\
&=|\iint(|\nabla \psi_R| \psi_R f^2)  drd\theta ds|\\
&=|\iint[  \partial_r (|\nabla \psi_R| \psi_R) f^2 +  |\nabla
\psi_R| \psi_R \partial_r (f^2) ] r \, drd\theta ds|.
\end{aligned}
\end{equation}
Using  Young's inequality, we deduce
\begin{eqnarray}
\label{b3=1} \nonumber
&& | \frac{1}{2}\iint(b_3 \cdot\nabla f^2)\psi_R^2dyds|\\
&& \lesssim \frac{1}{(\sigma_1 -
  \sigma_2)^2R^2}\|f\|_{L^2_tL^2\big(P(\sigma_1R)\big)}^2 +
  \frac{1}{8}\iint_{P(\sigma_1R)}|\nabla(\psi_Rf)|^2dyds.
\end{eqnarray}

Plugging in the above three estimates (\ref{b1=1}), (\ref{b2=2})
and (\ref{b3=1}) on terms involving $b_i$, $i=1, 2, 3$ into
\eqref{E2-3}, we arrive at
\begin{eqnarray}
&&\sup_{- \sigma_1^2R^2 \leq t \leq 0}\int_{B(\sigma_1R)}\psi_R^2
  f^2(\cdot, t)dy + \iint_{P(\sigma_1R)}
  |\nabla(f\psi_R)|^2dyds\\\nonumber
&&\lesssim \frac{1 + \dot{E}(\sigma_1R, \sigma_2R,
  b_1)^2}{(\sigma_1 - \sigma_2)^2R^2}\iint_{P(\sigma_1R)}f^2dyds
  +\frac{ \Vert B \Vert^2_{L^\infty_t BMO}R^{\frac{5}{3}}}{(\sigma_1 -
  \sigma_2)^2R^2} \| f\|^2_{L^3_tL^3\big(P(\sigma_1R)\big)}.
\end{eqnarray} By H\"older inequality, this implies
\begin{eqnarray}\label{E2-4}
&& \sup_{- \sigma_1^2R^2  \leq t \leq 0}\int_{B(\sigma_1
R)}\psi_R^2
  f^2(\cdot, t)dy + \iint_{P(\sigma_1R)}
  |\nabla(f\psi_R)|^2dyds\\\nonumber
&&\qquad \lesssim \frac{K(\| b\|_E) R^{\frac{5}{3}}}{(\sigma_1 -
  \sigma_2)^2R^2} \| f\|^2_{L^3_tL^3\big(P(\sigma_1R)\big)}.
\end{eqnarray}
Here and later in the section, as has been
mentioned in the introduction, $K=K(\cdot)$ is a one variable
function which may change from line to line, and $\| b\|_E$ is
defined in (\ref{Enorm}).

Our next step is to derive a mean value inequality based on
\eqref{E2-4} using Moser's iteration method. By H${\rm
\ddot{o}}$lder inequality and Sobolev imbedding theorem, one has
\begin{eqnarray}\nonumber
&&\iint_{P(\sigma_1R)}(\psi_R
  f)^{\frac{10}{3}}dyds \lesssim \int\Big(\|f\psi_R(\cdot, s)\|_{L^2(B(\sigma_1R))
  }^{\frac{4}{3}}\|\nabla(f\psi_R)\|_{L^2(B(\sigma_1R))}^2\Big)ds\\\nonumber
&&\lesssim \sup_{ - (\sigma_1R)^2 \leq s < 0}
  \|f\psi_R(\cdot, s)\|_{L^2(B(\sigma_1R))}^{\frac{4}{3}}
  \|\nabla(f\psi_R)\|_{L^2(P(\sigma_1R))}^2.
\end{eqnarray}
Using \eqref{E2-1} and \eqref{E2-4}, we obtain
\begin{eqnarray}\nonumber
\iint_{P(\sigma_2R)}f^{\frac{10}{3}}dyds \lesssim
  \Big\{\frac{K(\| b\|_E)^{\frac{3}{2}}}{(\sigma_1 - \sigma_2)^3R^{\frac{1}{2}}}
  \iint_{P(\sigma_1R)}f^3 dyds\Big\}^{\frac{10}{9}},
\end{eqnarray}
which implies that
\begin{eqnarray}\label{E2-5}
\iint_{P(\sigma_2R)}(|\Gamma|^{3q})^{\frac{10}{9}}dyds \leq
  \Big\{\frac{K(\| b\|_E)^{\frac{3}{2}}}{(\sigma_1 - \sigma_2)^3R^{\frac{1}{2}}}
  \iint_{P(\sigma_1R)}|\Gamma|^{3q}dyds\Big\}^{\frac{10}{9}}.
\end{eqnarray}

For integers $j \geq 0$ and a constant $\sigma = \frac{1}{3}$, set
$\sigma_2 = \frac{1}{2}\big(1 + \sigma^{j + 1}\big)$, $\sigma_1 =
\frac{1}{2}\big(1 + \sigma^{j}\big)$, $q = (\frac{10}{9})^j$ in
\eqref{E2-5}. Then we have
\begin{eqnarray}\nonumber
&&\Big\{\iint_{P\Big(\frac{R}{2}\big(1 + \frac{1}{\sigma^{j +
  1}}\big)\Big)}|\Gamma|^{3(\frac{10}{9})^{j + 1}}dyds
  \Big\}^{\frac{1}{3}(\frac{9}{10})^{j + 1}}\\\nonumber
&&\leq \Big\{\frac{K(\| b\|_E)^{\frac{3}{2}}}
  {\sigma^{3j}R^{\frac{1}{2}}}
  \iint_{P\Big(\frac{R}{2}\big(1 + \frac{1}{\sigma^{j}}\big)\Big)}
  |\Gamma|^{3(\frac{10}{9})^j}dyds\Big\}^{\frac{1}{3}(\frac{9}{10})^{j}}.
\end{eqnarray}
By iteration, the above inequality gives
\begin{eqnarray}\nonumber
&&\Big\{\iint_{P\Big(\frac{R}{2}\big(1 + \frac{1}{\sigma^{j +
  1}}\big)\Big)}|\Gamma|^{3(\frac{10}{9})^{j + 1}}dyds
  \Big\}^{\frac{1}{3}(\frac{9}{10})^{j + 1}}\\\nonumber
&&\leq \Big\{\frac{K(\| b\|_E)^{\frac{3}{2}}}
  {R^{\frac{1}{2}}}\Big\}^{\sum_{k =
  0}^j\frac{1}{3}(\frac{9}{10})^{k}}\sigma^{- \sum_{l =
  0}^jl(\frac{9}{10})^{l}}\Big\{\iint_{P(R)}|\Gamma|^3dyds\Big\}^{\frac{1}{3}}\\\nonumber
&&\lesssim \Big\{\frac{K(\|b\|_E)^{\frac{3}{2}}}{R^{\frac{1}{2}}}
  \Big\}^{\frac{10}{3}\big(1 - (\frac{9}{10})^{j}\big)}
  \Big\{\iint_{P(R)}|\Gamma|^3dyds\Big\}^{\frac{1}{3}}.
\end{eqnarray}
We take the limit $j \rightarrow \infty$ to yield that
\begin{eqnarray}\label{E2-6}
\sup_{P\big(\frac{R}{2}\big)}|\Gamma| \lesssim \big( K(\| b \|_E)
\big)^5\Big\{\frac{1}{R^5}\iint_{P(R)}
|\Gamma|^3dyds\Big\}^{\frac{1}{3}}.
\end{eqnarray}
From this a well known algebraic trick (see p87 \cite{HaL} e.g.)
shows
\begin{eqnarray}\label{E2-7}
\sup_{P\big(\frac{R}{2}\big)}|\Gamma| \lesssim \big( K( \| b \|_E)
\big)\Big\{\frac{1}{R^5}\iint_{P(R)}
|\Gamma|^2dyds\Big\}^{\frac{1}{2}}.
\end{eqnarray} Here the function $K(\cdot)$ may have changed at
the last step.

\section{H${\rm \ddot{o}}$lder Continuity of\ \ $\Gamma$}

In this section we study the regularity of $\Gamma$ using the
local maximum estimates of \eqref{E2-6} in section 2 and Nash type
method for parabolic equations.

Let us first recall a Nash inequality, whose proof can be found in
\cite{CSYT2009}.
\begin{lem}\label{Lem3-1}
Let $M \geq 1$ be a constant and $\mu$ be a probability measure.
Then for all $0 \leq f \leq M$, there holds
\begin{equation}\nonumber
\Big|\ln\int fd\mu - \int\ln fd\mu\Big| \leq
\frac{M\|g\|_{L^2}}{\int fd\mu},
\end{equation}
where $g = \ln f - \int\ln fd\mu$.
\end{lem}

Let $\zeta$ be a smooth radial cut-off function such that
\begin{equation}\label{3.1}
\begin{cases}
\zeta = 1\ \ {\rm on}\ \ B(\frac{1}{2}),\quad \zeta = 0\ \ {\rm
  on}\ \ B(1)^c,\quad \int_{\mathbb{R}^3}\zeta^2(x)dx = 1,\\[-4mm]\\
\big|\frac{\nabla\zeta}{\sqrt{\zeta}}\big| < \infty,\quad
\big|\nabla\big(\frac{\nabla\zeta}{\sqrt{\zeta}}\big)\big| <
\infty,
\end{cases}
\end{equation}
and $\zeta_R(x) = \frac{1}{R^{\frac{3}{2}}}\zeta(\frac{x}{R})$.
Let $\Phi$ be a positive solution to \eqref{1.5} in $P(R)$.

\begin{lem}\label{Lem3-2}
Let $\Phi \leq 2$ be a positive solution to \eqref{1.5} in $P(R)$
which is assumed to satisfy
\begin{equation}\label{E3-1}
\|\Phi\|_{L^1\big(P(\frac{R}{2})\big)} \geq c_0R^5.
\end{equation}
Moreover, we assume that $\Phi(r = 0, z, t)$ is a constant bigger
than 1. Then there holds
\begin{equation}\label{E3-2}
-\int\zeta_R^2(x)\ln\Phi(x, t)dx \leq M_0\big(1 + \|b\|_{E}^2\big)
\end{equation}
for all $t \in [-  \frac{c_0R^2}{4}, 0]$ and some absolute
positive constant $M_0$ depending only on $c_0$.
\end{lem}

\begin{proof}
First of all, let us define $\widetilde{\Phi}(x, t) = \Phi(Rx,
R^2t)$ and $\widetilde{b}(x, t) = Rb(Rx, R^2t)$. It is clear that
$\widetilde{\Phi}$ solves the equation
\begin{equation}\nonumber
\partial_t\widetilde{\Phi} +
\widetilde{b}\cdot\nabla\widetilde{\Phi} +
\frac{2}{r}\partial_r\widetilde{\Phi} = \Delta\widetilde{\Phi}
\end{equation}
on $P(1)$ and $0 \leq \widetilde{\Phi} \leq 2$,
$\|\widetilde{\Phi}\|_{L^1\big(P(\frac{1}{2})\big)} \geq c_0$. The
quantity we are going to control is $-\int\zeta_R^2(x)\ln\Phi(x,
t)dx = -\int\zeta^2(x)\ln\widetilde{\Phi}(x, R^{-2}t)dx$ on a time
internal $[-  \frac{c_0R^2}{4}, 0]$. Equivalently, we just need
to estimate $-\int\zeta^2(x)\ln\widetilde{\Phi}(x, t)dx$ for $t
\in [- \frac{c_0}{4}, 0]$.

Let $\Psi = - \ln\widetilde{\Phi}$. It is easy to see that $\Psi$
solves the equation
\begin{equation}\label{E3-3}
\partial_t\Psi + \widetilde{b}\cdot\nabla\Psi +
\frac{2}{r}\partial_r\Psi - \Delta\Psi + |\nabla\Psi|^2 = 0.
\end{equation}
Hence, by testing \eqref{E3-3} with $\zeta^2$ and using
integrating by parts and Cauchy-Schwarz's inequality, one has
\begin{eqnarray}\nonumber
&&\partial_t\int\Psi\zeta^2dx + \int|\nabla\Psi|^2\zeta^2dx  =
  \int\big(- \widetilde{b}\cdot\nabla\Psi - \frac{2}{r}\partial_r\Psi + \Delta\Psi\big)
  \zeta^2dx\\\nonumber
&&\leq - \int\zeta^2\widetilde{b}\cdot\nabla\big(\Psi -
  \bar{\Psi}(s)\big)dx - \int\frac{4\pi}{r}\partial_r\big(\Psi -
  \bar{\Psi}(s)\big)\zeta^2rdr d\theta dz\\\nonumber
&&\quad +\ \frac{1}{4}\int|\nabla\Psi|^2\zeta^2dx +
\int|\nabla\zeta|^2dx.
\end{eqnarray}
Here $\bar{\Psi}(s) = \int\Psi(\cdot, t)\zeta^2dx$. Using the
weighted Poincar${\rm \acute{e}}$ inequality
\begin{eqnarray}\label{E3-4}
\int|\Psi - \bar{\Psi}(s)|^2\zeta^2dx \leq
C\int|\nabla\Psi|^2\zeta^2dx
\end{eqnarray}
and the divergence-free property of $b$, we can estimate
\begin{eqnarray}\nonumber
\Big|\int\zeta^2\widetilde{b_1}\cdot\nabla\big(\Psi - \bar{\Psi}
(s)\big)dx\Big| \leq \frac{1}{8}\int|\nabla\Psi|^2\zeta^2dx +
C\int|\widetilde{b_1}|^2|\nabla\zeta|^2dx,
\end{eqnarray}
and
\begin{eqnarray}\nonumber
&&\Big|\int\zeta^2\widetilde{b_2}\cdot\nabla\big(\Psi -
  \bar{\Psi} (s)\big)dx\Big| =
  \Big|\int\nabla\zeta^2 \nabla\times B\big(\Psi -
  \bar{\Psi} (s)\big)dx\Big|\\\nonumber
&&\lesssim \int|\zeta\nabla(\Psi - \bar{\Psi})|
  |\nabla\zeta||B - \bar{B}|dx + \int\big|\nabla\big(\zeta^{\frac{3}{2}}\frac{\nabla\zeta}
  {\sqrt{\zeta}}\big)\big||B - \bar{B}||\Psi -  \bar{\Psi} (s)|dx\\\nonumber
&&\lesssim \int|\zeta\nabla(\Psi - \bar{\Psi})|
  |\nabla\zeta||B - \bar{B}|dx + \int\big(\frac{|\nabla\zeta|^2}{\zeta} + \sqrt{\zeta}\nabla
  \frac{\nabla\zeta}{\sqrt{\zeta}}\big)|B - \bar{B}||\zeta||\Psi - \bar{\Psi} (s)|dx\\\nonumber
&&\lesssim \frac{1}{8}\int|\nabla\Psi|^2\zeta^2dx + C\|B\|_{{\rm
  BMO}}^2.
\end{eqnarray} Here we just used the weighted Poincar\'e
inequality and (\ref{BMOLp}), with $p=2$. Moreover
\begin{eqnarray}\nonumber
&&\Big|\int\zeta^2\widetilde{b_3}\cdot\nabla\big(\Psi -
  \bar{\Psi} (s)\big)dx\Big| =
  \Big|\int\nabla\zeta^2\widetilde{b_3}\big(\Psi -
  \bar{\Psi} (s)\big)dx\Big|\\\nonumber
&&\lesssim \int|\nabla\zeta(\Psi - \bar{\Psi})|\zeta dr d\theta dz
  = \int|\frac{\nabla\zeta}{\sqrt{\zeta}}(\Psi - \bar{\Psi})|\zeta^{\frac{3}{2}} dr d\theta dz\\\nonumber
&&\lesssim \int|\zeta\nabla(\Psi - \bar{\Psi})|
  |\nabla\zeta|dx + \int\big(\frac{|\nabla\zeta|^2}{\zeta} + \sqrt{\zeta}\nabla
  \frac{\nabla\zeta}{\sqrt{\zeta}}\big)|\zeta||\Psi - \bar{\Psi} (s)|dx\\\nonumber
&&\lesssim \frac{1}{8}\int|\nabla\Psi|^2\zeta^2dx + C.
\end{eqnarray}
Here we also used the integration by parts. On the other hand, by
recalling the assumption that $\Phi(r = 0, z, t)$ is a non-zero
constant, one can estimate
\begin{eqnarray}\nonumber
&&- \int\frac{4\pi}{r}\partial_r\Psi\zeta^2rdr d\theta dz =
  - 4\pi\int_{-\infty}^\infty(\Psi - \bar{\Psi})\zeta^2dz\Big|_{r = 0}^{r = \infty}
  + 4\pi\int(\Psi - \bar{\Psi})\partial_r\zeta^2dr d\theta dz\\\nonumber
&&= 4\pi\int_{-\infty}^\infty\Psi\zeta^2 dz\Big|_{r = 0} -
  4\pi\bar{\Psi}\int_{-\infty}^\infty\zeta^2dz
  \Big|_{r = 0} + 4\pi\int(\Psi - \bar{\Psi})\zeta\frac{\partial_r\zeta}{r}rdr d\theta dz\\\nonumber
&&\leq C - C\bar{\Psi}(s) +
\frac{1}{8}\int|\nabla\Psi|^2\zeta^2dx.
\end{eqnarray}
Here we also used the fact that the support of
$\frac{1}{r}|\partial_r\zeta|$  is away from $z$- axis.
Consequently, we obtain
\begin{eqnarray}\nonumber
\partial_t\int\Psi\zeta^2dx + C\int\Psi\zeta^2dx \leq -
\frac{1}{2}\int|\nabla\Psi|^2\zeta^2dx + C\big(1 +
\|b\|_{E}^2\big).
\end{eqnarray}

In order to proceed, we apply the Nash inequality in Lemma
\ref{Lem3-1}. Take $f = \widetilde{\Phi}$, $d\mu = \zeta^2(x)dx$.
One has
\begin{equation}\nonumber
\Big|\ln\int\widetilde{\Phi}\zeta^2dx +
\int\Psi\zeta^2dx\Big|^2\Big(\int \widetilde{\Phi}\zeta^2dx\Big)^2
\leq M^2\int\big|- \Psi + \int\Psi\zeta^2dy\big|^2\zeta^2dx.
\end{equation}
Here $M = 2$ is the upper bound of $\Phi$. Using the weighted
Poincar${\rm \acute{e}}$ inequality \eqref{E3-4} once again, we
have
\begin{equation}\nonumber
\Big|\ln\int\widetilde{\Phi}\zeta^2dx +
\int\Psi\zeta^2dx\Big|^2\Big(\int \widetilde{\Phi}\zeta^2dx\Big)^2
\leq C\int|\nabla\Psi|^2\zeta^2dx.
\end{equation}
Hence, we finally obtain
\begin{eqnarray}\nonumber
&&\partial_t\bar{\Psi}(t) + C_0\bar{\Psi}(t) \leq
  C\big(1 + \|b\|_{E}^2\big)\\\nonumber &&\quad -\
(2C)^{-1}\Big|\ln\int\widetilde{\Phi}\zeta^2dx +
  \bar{\Psi}\Big|^2\Big(\int \widetilde{\Phi}\zeta^2dx\Big)^2.
\end{eqnarray}
Let $\chi(s)$ be the characteristic function of the set
\begin{eqnarray}\nonumber
W = \big\{s \in [-\frac{1}{4}, 0):
\|\widetilde{\Phi}(s)\|_{L^1\big(B_{\frac{1}{2}}\big)} \geq
\frac{c_0}{2}\big\}.
\end{eqnarray}
By the assumption \eqref{E3-1} and hence
$\|\widetilde{\Phi}\|_{L^1\big(P(\frac{1}{2})\big)} \geq c_0$, one
has $|W| \geq \frac{3c_0}{4}$. In fact, if $|W| < \frac{3c_0}{4}$,
then
\begin{equation}\nonumber
\|\widetilde{\Phi}\|_{L^1\big(P(\frac{1}{2})\big)} <
\int_{W}2\big|B\big(\frac{1}{2}\big)\big|ds +
\int_{W^c}\frac{c_0}{2}ds \leq \frac{(2\pi + 1)c_0}{8} < c_0,
\end{equation}
which contradicts with \eqref{E3-1}. Thus, we have
\begin{eqnarray}\label{E3-5}
\partial_t\bar{\Psi}(t)  + C_0\bar{\Psi}(t) \leq
  C_0\big(1 + \|b\|_{E}^2\big) - 8C_0^{-1}c_0^2\chi(s)\Big|\ln\int
  \widetilde{\Phi}\zeta^2dx + \bar{\Psi}\Big|^2.
\end{eqnarray}
Note that the obvious consequence of this inequality gives
\begin{eqnarray}
\label{psis1s2}
 \bar{\Psi}(s_2) \leq \bar{\Psi}(s_1) +
C_0e^{C_0}(s_2 - s_1)\big(1 + \|b\|_{E}^2\big)
\end{eqnarray}
for $- \frac{1}{4} \leq s_1 \leq s_2 \leq 0$. Hence, if for some
$s_0 \in [- \frac{1}{4}, - \frac{c_0}{4})$ such that
\begin{eqnarray}\nonumber
\bar{\Psi}(s_0) \leq \frac{4C_0}{c_0}\big(1 + \|b\|_{E}\big) +
2\big|\ln\frac{c_0}{2}\big|,
\end{eqnarray}
then we are done since
\begin{eqnarray}\nonumber
\bar{\Psi}(t) \leq \bar{\Psi}(s_0) + \frac{C_0e^{C_0}}{2}\big(1 +
\|b\|_{E}^2\big)
\end{eqnarray}
for all $t \in [s_0, 0)$. Otherwise, one has
\begin{eqnarray}\nonumber
\bar{\Psi}(s) \geq \frac{4C_0}{c_0}\big(1 + \|b\|_{E}\big) +
2\big|\ln\frac{c_0}{2}\big|
\end{eqnarray}
for all $s \in [- \frac{1}{4} , - \frac{c_0}{4})$. For $s \in W
\cap [- \frac{1}{4} , - \frac{c_0}{4})$, one has
\begin{eqnarray}\nonumber
\ln\int\widetilde{\Phi}\zeta^2dx \geq
\ln\int_{B_{\frac{1}{2}}}\widetilde{\Phi}dx \geq \ln\frac{c_0}{2}.
\end{eqnarray}
Hence, by \eqref{E3-5}, we have
\begin{eqnarray}\nonumber
\partial_t\bar{\Psi} + C_0\bar{\Psi} \leq - \frac{c_0^2\chi(s)}{C_0}\bar{\Psi}^2,\quad
- \frac{1}{4} \leq s \leq - \frac{c_0}{4}.
\end{eqnarray}
Solving the above inequality gives
\begin{eqnarray}\nonumber
\bar{\Psi}(- \frac{c_0}{4}) \leq
\frac{1}{\frac{c_0^2}{C_0}\int_{-\frac{1}{4}}^{-\frac{c_0}{4}}\chi(s)e^{-C_0s}ds
+ \frac{1}{\overline{\Psi}(- \frac{1}{4})}} < \infty.
\end{eqnarray}
The bound of the $\bar{\Psi}(- \frac{c_0}{4})$ depends only on
$c_0$ since $\overline{\Psi}(- \frac{1}{4}) > 0$ and $|W| \geq
\frac{3c_0}{4}$. Starting from $s = - \frac{c_0}{4}$ and using
(\ref{psis1s2}), we have
\begin{eqnarray}\nonumber
\bar{\Psi}(s) \leq \bar{\Psi}\big( - \frac{c_0}{4}\big) +
C_0e^{C_0}\big(1 + \|b\|_{E}^2\big)
\end{eqnarray}
for all $s \in [- \frac{c_0}{4}, 0]$, which completes the proof of
the lemma.
\end{proof}

As an corollary, Lemma \ref{Lem3-2} gives a lower bound of
positive solutions of \eqref{1.5}.
\begin{cor}\label{Cor3-1}
Let $\Phi$, $c_0$ and $M_0$ be given in Lemma \ref{Lem3-2}. Then
there exists a constant $0 < \delta < 1$ depending only on
$\|b\|_{E}$ such that
\begin{equation}\label{E3-6}
\inf_{P\big(\frac{R}{8}\big)}\Phi \geq \frac{\delta}{2}.
\end{equation}
\end{cor}
\begin{proof}
Using Lemma \ref{Lem3-2}, we have
\begin{eqnarray}\nonumber
&&M_0\big(1 + \|b\|_{E}^2\big) \geq
  - \int\zeta_R^2(x)\ln\Phi(t, x)dx\\\nonumber
&&= -\int_{\delta < \Phi \leq 1}\zeta_R^2(x)\ln\Phi(t, x)dx -
  \int_{\Phi \leq \delta}\zeta_R^2(x)\ln\Phi(t, x)dx\\\nonumber
&&\quad -\ \int_{1 < \Phi \leq 2}\zeta_R^2(x)\ln\Phi(t,
  x)dx\\\nonumber
&&\geq 0 - \int_{\Phi \leq \delta}\zeta_R^2(x)\ln\Phi(t, x)dx -
  \ln2\int_{1 < \Phi \leq 2}\zeta_R^2(x)dx\\\nonumber
&&\geq - \int_{\Phi \leq \delta}\zeta_R^2(x)\ln\Phi(t, x)dx  -
  \ln2,
\end{eqnarray}
which implies that
\begin{eqnarray}\nonumber
&&- \int_{\Phi \leq \delta}\zeta_R^2(x)\ln\Phi(t, x)dx \lesssim 1
+ \|b\|_{E}^2
\end{eqnarray}
for $- \frac{c_0 R^2}{4} \leq t \leq 0$. Consequently, we have
\begin{eqnarray}\nonumber
\big|\big\{x \in B\big(\frac{R}{2}\big)| \Phi(t, x) \lesssim
\delta\big\}\big| \leq \frac{R^3}{- \ln\delta}\big(1 +
\|b\|_{E}^2\big)
\end{eqnarray}
for $- \frac{R^2}{64} \leq t \leq 0$. Using the mean value
inequality \eqref{E2-7}, one has
\begin{eqnarray}\nonumber
\sup_{P\big(\frac{R}{8}\big)}(\delta - \Phi)_+
  &\lesssim& \Big\{\frac{K(\|b\|_E)}{R^5}\iint_{P(\frac{R}{2})}
  (\delta - \Phi)_+^2dyds\Big\}^{\frac{1}{2}}\\\nonumber
&\lesssim& \frac{\delta}{\sqrt{|\ln\delta|}}K(\|b\|_E),
\end{eqnarray}
which gives
\begin{eqnarray}\nonumber
\inf_{P\big(\frac{\sqrt{c_0}R}{2}\big)}\Phi \geq \delta -
\frac{C_0\delta}{2\sqrt{|\ln\delta|}} K(\|b\|_E)
\end{eqnarray}
for some $C_0 > 0$ which is independent of $\delta$ and $R$. Then
\eqref{E3-6} follows by choosing a sufficiently small $\delta$
such that
\begin{eqnarray}\label{E3-7}
\delta \leq \exp\big\{- K(\|b\|_{E})\big\}.
\end{eqnarray}
\end{proof}

Now we are ready to give the proof of Theorem \ref{Holder}.

{\bf Proof of Theorem \ref{Holder}}.

Without loss of generality we take $L=1$.

For $0 < r \leq 1$, we define
\begin{eqnarray}\nonumber
m_r = \inf_{P(r)}\Gamma,\quad M_r = \sup_{P(r)}\Gamma,\quad J_r =
M_r - m_r.
\end{eqnarray}
As in \cite{CSYT2009}, we define
\begin{eqnarray}\nonumber
\Phi = \begin{cases}\frac{2(M_1 - \Gamma)}{J_1}\quad {\rm if}\ M_1 > - m_1,\\
\frac{2(\Gamma - m_1)}{J_1}\quad {\rm if}\ M_1 \leq -
m_1.\end{cases}
\end{eqnarray}
It is clear that $0 \leq \Phi \leq 2$ is a non-negative solution
of \eqref{1.5} in $P(1)$ and $a = \Phi(r = 0, z, t)$ is a constant
bigger than 1. To verify that $\Phi$ satisfies the condition
\eqref{E3-1}, we need the following lemma on the lower bound of
$\|\Phi\|_{L^p}$ for $0 < p < 1$ as in \cite{CSYT2009}.

\begin{lem}\label{lemLB}
Suppose that $b$ satisfies \eqref{1.1} and \eqref{1.2}. Then for
arbitrary $p \in (0, 1)$, $\Phi$ defined above satisfies
\begin{equation}\nonumber
\frac{1}{R^\frac{5}{p}}\|\Phi\|_{L^{p}(P(R, \frac{R}{2}))} \geq
C^{-1}\big( K(\|b\|_{E})\big)^{- \frac{2}{p}}a.
\end{equation}
\end{lem}
\begin{proof}

Since the lemma is scaling invariant, we just take $R=1$ in the
proof.
 Let $\psi = \phi(|x|)\eta(t)$, where $\phi \in C_0^\infty$
such that $\phi = 1$ on $B_{\frac{1}{2}}$, $\phi = 0$ on $B_1^c$,
$\frac{\nabla\phi}{\sqrt{\phi}}$ and
$\nabla\frac{\nabla\phi}{\sqrt{\phi}}$ are bounded, $\eta \in
C_0^\infty$ such that $\eta = 1$ on $[- \frac{7}{8}, -
\frac{1}{8}]$ and $\eta $ is supported in $(- 1, 0)$. Let us test
\eqref{1.5} by $p\Phi^{p - 1}\psi_R^2$, $p \in (0, \frac{1}{2})$,
to derive that
\begin{eqnarray}
\label{lemma3.4-1}
 \iint\Big(\partial_s\Phi^p +
(b\cdot\nabla)\Phi^p + \frac{2}{r}\partial_r\Phi^p\Big)\psi^2dyds
= p\iint
  \Delta\Phi\Phi^{p - 1}\psi^2dyds.
\end{eqnarray}
Similarly as in \cite{CSYT2009}, we have
\begin{eqnarray}
\label{2/rsec3} &&-
\iint\frac{2}{r}(\partial_r\Phi^p)\psi^2dyds\\\nonumber &&=
  \iint\Phi^p\frac{4}{|y'|}\psi(\partial_{|y'|}\psi)dyds +
  \int_{-1}^0ds\int2\Phi^p\psi^2\big|_{r = 0}dz\\\nonumber
&&\geq - C \iint\Phi^pdyds + \frac{3}{2}a^p.
\end{eqnarray} Here $|y'|=\sqrt{y^2_1+y^2_2}$ if $y=(y_1, y_2,
y_3)$. Likewise
\begin{eqnarray}\label{ddphisec3}
&&\iint\big(- \partial_s\Phi^p  + p\Delta\Phi\Phi^{p -
  1}\big)\psi^2dyds\\\nonumber
&&= \iint2\Phi^p\big[\psi(\partial_s\psi) + |\nabla\zeta|^2 -
  \frac{p - 2}{p}\zeta\Delta\zeta\big]dyds - \frac{4(p -
  1)}{p}\iint|\nabla(\Phi^{\frac{p}{2}}\psi)|^2dyds\\\nonumber
&&\geq - C \iint\Phi^pdyds - \frac{4(p -
  1)}{p}\iint|\nabla(\Phi^{\frac{p}{2}}\psi)|^2dyds.
\end{eqnarray}
Moreover, concerning the term involving $b$, we estimate it as
follows:
\begin{eqnarray}
\label{b1sec3} &&- \iint\psi^2(b_1\cdot\nabla)\Phi^pdyds =
  \iint\Phi^pb_1\cdot\nabla\psi^2dyds\\\nonumber
&&\geq - C \|b_1\|_{L^\infty_tL^{2}(P(1, \frac{1}{2}))}
  \|\Phi^p\|_{L^1_tL^{2}(P(1, \frac{1}{2}))},
\end{eqnarray}
and
\begin{eqnarray}\nonumber
&&- \iint\psi^2(b_2\cdot\nabla)\Phi^pdyds =
  \iint (B-\bar{B}) \cdot\nabla\times(\Phi^p\nabla\psi^2)dyds\\\nonumber
&&= \iint (B-\bar{B}) \cdot\nabla\times[(\Phi^{\frac{p}{2}}\psi)^2
  \frac{\nabla\psi}{\psi}]dyds\\\nonumber &&\geq -
  \frac{C}{R}\|\Phi^{\frac{p}{2}}\psi\|_{L^2_tL^{2}(P(1,
  \frac{1}{2}))} \|\nabla(\Phi^{\frac{p}{2}}\psi)\|_{L^2_tL^{2}(P(1,
  \frac{1}{2}))} - C \iint |B-\bar{B}| \Phi^pdyds\\\nonumber
&&\geq \frac{2(p - 1)}{p} \iint|\nabla(\Phi^{\frac{p}
  {2}}\psi)|^2dyds - C \iint\Phi^p dyds\\\nonumber
&&\quad -\  C \big{(} \iint\Phi^{2p} dyds \big{)}^{1/2}
  \big{(} \iint (B-\bar{B})^2 dyds \big{)}^{1/2}.
\end{eqnarray}
By H\"older inequality and (\ref{BMOLp}), we have
\begin{equation}
\label{b2sec3} - \iint\psi^2(b_2\cdot\nabla)\Phi^pdyds \ge
\frac{2(p - 1)}{p} \iint|\nabla(\Phi^{\frac{p}
  {2}}\psi)|^2dyds - K(\|b\|_E) \|\Phi\|_{L^{2p}(P(1,
  \frac{1}{2}))}^p
\end{equation}
Just like (\ref{b3=1}), we also have

\begin{equation}
\label{b3sec3} -\iint\psi^2(b_3 \cdot\nabla)\Phi^pdyds
   \ge \frac{(p - 1)}{p} \iint|\nabla(\Phi^{\frac{p}
  {2}}\psi)|^2dyds - C \iint\Phi^pdyds.
\end{equation}

Substituting (\ref{b3sec3}), (\ref{b2sec3}), (\ref{b1sec3}),
(\ref{ddphisec3}) and  (\ref{2/rsec3}) into (\ref{lemma3.4-1}), we
deduce
\begin{eqnarray}\nonumber
&&\frac{3}{2}a^p \leq C \iint\Phi^pdyds + C
 K(\| b\|_E)
  \|\Phi^p\|^p_{L^2_tL^{2}(P(1, \frac{1}{2}))}\\\nonumber
&&\leq C K(\| b\|_E) \|\Phi\|_{L^{2p}(P(1,
  \frac{1}{2}))}^p,
\end{eqnarray}
which completes the proof of the lemma, since $p \in (0, 1/2)$ is
arbitrary.
\end{proof}

Now we continue the proof of the theorem.

By Lemma \ref{lemLB}, $\Phi$ satisfies the assumptions in Lemma
\ref{Lem3-2} for $R = 1$. By Corollary \ref{Cor3-1}, one has
\begin{equation}\nonumber
\inf_{P\big(\frac{\sqrt{c_0}}{2}\big)}\Phi \geq \frac{\delta}{2}.
\end{equation}
Noting that and $m_1 \leq
\inf_{P\big(\frac{\sqrt{c_0}}{2}\big)}\Gamma \leq
\sup_{P\big(\frac{\sqrt{c_0}}{2}\big)}\Gamma \leq M_1$, we have
\begin{equation}\label{E3-8}
J_{\frac{\sqrt{c_0}}{2}} = {\rm
OSC}_{P\big(\frac{\sqrt{c_0}}{2}\big)}\Gamma \leq \big(1 -
\frac{\delta}{4}\big)J_1.
\end{equation}
Iterating \eqref{E3-8} immediately shows that $\Phi$ is H${\rm
\ddot{o}}$lder continuous at $(0, 0)$. \qed

\section{Applications to axially symmetric Navier-Stokes equation}

This section is devoted to proving Theorems \ref{Liouville} and
\ref{thregcon}. We begin with
\medskip

\noindent{\bf Proof of Theorem \ref{Liouville}.}

By the assumptions of the theorem, we can apply Theorem
\ref{Holder} to deduce that the function $\Gamma = r v^{\theta}$
is H\"older continuous at the space time point $(0, 0)$. More
precisely, for any fixed point $(x, t) \in {\bf R}^3 \times
(-\infty, 0)$, there exist positive constants $\alpha$ and $C$
such that for all sufficiently large $L>0$, we have
\[
|\Gamma(x,t) -\Gamma(0, 0)| \le C [(|x|+\sqrt{|t|})/L]^{\alpha}
\sup \Gamma.
\]Letting $L \to \infty$, we find that
\[
r v^{\theta}(x, t) = \Gamma(x, t) = \Gamma(0, 0).
\]Since $v^{\theta}$ is a bounded function, the only way this can
happen is $v^{\theta} \equiv 0$. Hence $v$ is a bounded, weak
ancient solution without swirl. According to Theorem 5.2 in
\cite{KNSS}, the ancient solution $v=(0, 0, l(t))$ where $l=l(t)$
depends only on time. Therefore its stream function $B$ is a
harmonic function since $\Delta B = - \nabla \times v=0$. Since
the function $B=B(\cdot, t)$ is BMO, by (\ref{BMOLp}) we know
\[
\int_{|x|<R} |B(x, t)-\bar{B}(t)| dx \le C \|B(\cdot, t)\|_{BMO}
R^3.
\]Here $\bar{B}(t)$ is the average of $B(\cdot, t)$ in the ball
$B_R$. Since $B(\cdot, t)$ is harmonic, the mean value theorem
tells us that $\bar{B}(t)=B(0, t)$. Hence
\[
\int_{|x|<R} |B(x, t)| dx \le C \|B(\cdot, t)\|_{BMO} R^3 + B(0,
t) R^3.
\]The mean value theorem then implies that $B(\cdot, t)$ is a
bounded function since
\begin{eqnarray}\nonumber
&&|B(y, t)| = \Big|\frac{1}{4\pi|y|^3}\int_{B(y, |y|)}B(z,
  t)dz\Big|\\\nonumber
&&\leq \frac{1}{4\pi|y|^3}\int_{B(0, 2|y|)}|B(z, t)|dz \lesssim
  \|B(\cdot, t)\|_{BMO} + |B(0, t)|.
\end{eqnarray}

The classical Liouville theorem shows that the stream function
$B$, being a bounded function, is constant. Therefore $v= \nabla
\times B=0$. \qed

\medskip

\noindent{\bf Proof of Theorem \ref{thregcon}.} We use the method
of contradiction.  If there is a singularity to the
axially-symmetric Navier-Stokes equations \eqref{axi-NS}, then we
can generate a nonzero, bounded, weak ancient solution as in
\cite{KNSS}. Our Theorem \ref{Holder} and a scaling argument will
then be used to show that such a bounded ancient solution is
identically zero. This contradiction proves that singularity can
not occur.

By time shifting, we assume that the solution $v$ exists in the
time interval $[-1, 0]$ and that $t = 0$ is a blow up time of $v$.
The partial regularity theory in \cite{CKN, Lin} says that the
Hausdorff measure of the singular space-time set of any suitable
weak solution is zero. This implies that for axially-symmetric
Navier-Stokes equations \eqref{axi-NS}, suitable weak solutions
can only develop singularities on the symmetric axis $r = 0$.
Hence, without loss of generality, we may assume that $(0, 0)$ is
the earliest blowup point.

For $k \geq 1$, let $(x_k, t_k)$ be a sequence of points such that
\begin{equation}\label{1.4}
- 1 < t_k \nearrow 0,\quad Q_k = |v(x_k, t_k)| = \gamma_k\max_{-1
< t < t_k}|v(x, t)| \nearrow \infty,\quad \gamma_k \longrightarrow
1.
\end{equation}
Define a sequence of functions $\{v^{(k)}\}$ by
\begin{equation}\label{1.6}
v^{(k)}(x, t) =\frac{1}{Q_k}v(x_k + \frac{x}{Q_k}, t_k +
\frac{t}{Q_k^2}),\quad - Q_k^2(1 + t_k) \leq t \leq 0.
\end{equation}

 By \cite{BZ}, one can assume that $r_k = r(x_k)$ are uniformly
bounded. It is clear that $\{v^{(k)}\}$ defined in \eqref{1.6} are
mild solutions. Moreover, $\{v^{(k)}\}$ (up to a subsequence)
converges to a bounded ancient weak solution $u(x, t)$ to the
Navier-Stokes equations (for details, see \cite{KNSS}). By the
construction, $|u(x, t)| \leq 1$ and $|u(0, 0)| = 1$.

We consider two cases.

Case 1 is  when $r_k|v(x_k, t_k)| = r_kQ_k$ are uniformly bounded
by some positive constant $C$. Then the functions $\{v^{(k)}\}$
are also axi-symmetric with respect to an axis which is parallel
to the z-axis and is at distance at most $C$ from it.
Consequently, $u$ is also axi-symmetric with respect to a suitable
axis. Note that both the stream function and $r v^\theta$ are
scaling invariant. Thus the stream function of $u$ is in BMO and
$r u^\theta$ is also bounded.  Therefore we can apply Theorem
\ref{Holder} on $u$, which says that the swirl component of $u$
vanishes. By Theorem 5.2 in \cite{KNSS}, we conclude $u=(0, 0,
l(t))$ with $l=l(t)$ being a function of time only. But this show
$u=0$ as in the proof of the previous theorem. This contradiction
shows that Case 1 can not happen.

Case 2 is  when $r_k|v(x_k, t_k)| = r_kQ_k$ is not uniformly
bounded.

Hence, $r_kQ_k$ (up to a subsequence) goes to infinity as $k$
tends to infinity. Due to Caffarelli-Kohn-Nirenberg's partial
regularity theory, $\{x_k\}$ (up to a subsequence) converges to
$x_\ast$ which is a point on the z-axis such that $r_\ast = 0$.
Due to the axis symmetry of $v$, $x_k$ can be chosen so that
$\theta(x_k) \rightarrow \theta_\infty$ for a $\theta_\infty$.
Hence, $e_r(x_k) \rightarrow \nu$ and $e_\theta(x_k) \rightarrow
\nu^\perp=(-\nu_2, \nu_1, 0)$ for an unit vector $\nu = (\nu_1,
\nu_2, 0)$. Here $e_r(x)$ and $e_\theta(x)$ are defined as in the
introduction and $(r(x), \theta(x))$ is the polar coordinate of
$(x_1, x_2)$.

It is clear that
\begin{equation}\nonumber
\begin{cases}
x_k + \frac{x}{Q_k} \in B(x_k, \frac{r_k}{\sqrt{r_kQ_k}})\quad
{\rm for}\ x \in B(0, \sqrt{Q_kr_k}),\\[-4mm]\\
t_k - (\frac{r_k}{\sqrt{r_kQ_k}})^2 < t_k + \frac{t}{Q_k^2} \leq
t_k < 0\quad {\rm for}\quad - Q_kr_k  < t \leq 0.
\end{cases}
\end{equation}
By the assumption on initial value and the maximum principle, we
know
\begin{equation}\nonumber
|v^\theta(t, y)| \lesssim \frac{1}{r_k} \quad {\rm for}\ y \in
B(x_k, \frac{r_k}{2}),\ t < 0,
\end{equation}
which shows
\begin{equation}\label{3-4}
\begin{cases}
\big|v^{(k)}(x, t)e_\theta\big(x_k + \frac{x}{Q_k}\big)\big| =
  \frac{1}{Q_k}\big|v^\theta\big(x_k + \frac{x}{Q_k}, t_k +
  \frac{t}{Q_k^2}\big)\big| \lesssim \frac{1}{Q_kr_k},\\[-4mm]\\
\quad{\rm for}\ (x, t) \in B(0, \sqrt{r_kQ_k})
  \times (- r_kQ_k, 0].
\end{cases}
\end{equation}

Note that  on $B(0, \sqrt{Q_kr_k}) \times (- Q_kr_k, 0]$, it is
easy to see that $e_r\big(x_k + \frac{x}{Q_k}\big) \rightarrow
\nu$ and $e_\theta\big(x_k + \frac{x}{Q_k}\big) \rightarrow
\nu^\perp$ as $n \rightarrow \infty$. Moreover,  for each $k$,
$v^{(k)}$ is still a mild solution to the 3D Navier-Stokes
equations. By \eqref{3-4}, there exists a subsequence of
$\{v^{(k)}\}$ (we will still denote it by $\{v^{(k)}\}$) and a
bounded ancient solution $u(x, t)$ to the 3D Navier-Stokes
equations on $\mathbb{R}^3 \times (- \infty, 0]$, which is mild in
the sense of \cite{KNSS}, such that
\begin{eqnarray}\nonumber
&&v^{(k)}(x, t) = \frac{1}{Q_k}v^r\big(x_k + \frac{x}{Q_k}, t_k +
  \frac{t}{Q_k^2}\big)e_r\big(x_k+ \frac{x}{Q_k}\big)\\\nonumber
&&\quad +\ \frac{1}{Q_k}v^\theta\big(x_k + \frac{x}{Q_k}, t_k +
  \frac{t}{Q_k^2}\big)e_\theta\big(x_k + \frac{x}{Q_k}\big)\\\nonumber
&&\quad +\ \frac{1}{Q_k}v^z\big(x_k + \frac{x}{Q_k}, t_k +
  \frac{t}{Q_k^2}\big)e_z\big(x_k + \frac{x}{Q_k}\big)\\\nonumber
&&\quad \rightarrow u = u^r\nu + u^\theta \nu^\perp + u^ze_z\quad
  {\rm in}\ L^\infty(\Omega)
\end{eqnarray}
for any compact subset $\Omega$ of $\mathbb{R}^3 \times
\mathbb{R}$ and $u(x, t)\cdot\nu^\perp = 0$. Hence,
\begin{eqnarray}\label{3-5}
u(x, t) = u^r(x, t)\nu + u^z(x, t)e_z.
\end{eqnarray}

On the other hand, for $(y, s) \in B(x_k,
\frac{r_k}{\sqrt{r_kQ_k}}) \times \big[t_k -
(\frac{r_k}{\sqrt{r_kQ_k}})^2, t_k\big]$, one has
\begin{eqnarray}\nonumber
&&- \frac{1}{Q_k}[v^r(y, s)e_\theta(y) - v^\theta(y,
  s)e_r(y)]\\\nonumber
&&= \frac{1}{Q_k}\partial_\theta[v^r(y, s)e_r(y) + v^\theta(y,
  s)e_\theta(y) + v^z(y, s)e_z(y)]\\\nonumber
&&= \partial_\theta\big[v^{(k)}\big(Q_k(y - x_k), Q_k^2(s -
  t_k)\big)\big]\\\nonumber
&&= Q_k(\partial_\theta y\cdot\nabla) v^{(k)}\big(Q_k(y - x_k),
  Q_k^2(s - t_k)\big)\\\nonumber
&&= Q_k|y|\big(e_\theta(y)\cdot\nabla\big) v^{(k)}\big(Q_k(y -
  x_k),  Q_k^2(s - t_k)\big),
\end{eqnarray}
which gives that
\begin{eqnarray}
\label{vdu}
 &&\frac{1}{Q_k}\big[- v^r\big(x_k + \frac{x}{Q_k}, t_k
+
  \frac{t}{Q_k^2}\big)e_\theta\big(x_k + \frac{x}{Q_k}\big)\\\nonumber
&&\quad + v^\theta\big(x_k + \frac{x}{Q_k}, t_k +
  \frac{t}{Q_k^2}\big)e_r\big(x_k + \frac{x}{Q_k}\big)\big]\\\nonumber
&&\qquad = Q_k\big|x_k + \frac{x}{Q_k}\big|\Big(e_\theta\big(x_k +
  \frac{x}{Q_k}\big)\cdot\nabla\Big)v^{(k)}(x, t)
\end{eqnarray}
for $(x, t) \in B(0, \sqrt{Q_kr_k}) \times (- Q_kr_k, 0]$. Since
$r_k Q_k \to \infty$, we know $Q_k\big|x_k + \frac{x}{Q_k}\big|
\to \infty$ for fixed $x$. But the left hand side of (\ref{vdu})
is bounded by definition of $Q_k$. Hence, let $k \to \infty$, we
have
\begin{eqnarray}\label{3-6}
(\nu^\perp\cdot\nabla) u(x, t) = 0.
\end{eqnarray}
Note that the Navier-Stokes equations are invariant under
rotation. Without loss of generality, we set $\nu = e_1$ and
$\nu^\perp = e_2$. Consequently, the limit function
\begin{eqnarray}\nonumber
u(x, t) = u^r(x_1, z, t)e_1 + u^z(x_1, z, t)e_z,
\end{eqnarray}
is a bounded ancient solution to the 2D Navier-Stokes equations.
By Theorem 5.1 in \cite{KNSS}, the limit
\begin{eqnarray}\label{3-7}
u(x, t) = u^r(t)\nu + u^z(t)e_z
\end{eqnarray}
depends only on $t$ and  that $|u(0, 0)| = 1$. By the argument in
the proof of the previous theorem, the boundedness of the stream
function of $u$ in $BMO$ norm implies that $u=0$. This
contradiction shows that Case 2 can not occur either. Therefore
the assumption that $v$ becomes singular at $(0, 0)$ is false,
proving Theorem \ref{thregcon}. \qed

\section*{Acknowledgement}
The authors would like to thank Professor Fang-hua Lin for some
helpful discussions. The work was in part supported by NSFC
(grants No. 10801029 and 10911120384), FANEDD, Shanghai Rising
Star Program (10QA1400300), SGST 09DZ2272900 and SRF for ROCS,
SEM.


\end{document}